\documentclass[12pt]{amsart}

\usepackage{accents}
\usepackage{appendix}
\usepackage{amsfonts}
\usepackage{amsmath}
\usepackage{amssymb}	
\usepackage{amsthm,bm}
\usepackage{array,booktabs,multirow}
\usepackage{braket}
\usepackage{centernot}
\usepackage{cite}
\usepackage{comment}
\usepackage{dsfont}
\usepackage{mathalfa}
\usepackage[shortlabels]{enumitem}
\usepackage{etoolbox}
\usepackage{float}
\usepackage[hang, flushmargin]{footmisc}
\usepackage{latexsym}
\usepackage{lipsum}
\usepackage{needspace}
\usepackage{tikz}
\usepackage{hyperref}
\usetikzlibrary{matrix,arrows}

\makeatletter
\def\namedlabel#1#2{\begingroup
   \def\@currentlabel{#2}%
   \label{#1}\endgroup
}
\makeatother

\newcommand{\floor}[1]{\left\lfloor #1 \right\rfloor}
\newcommand{\ceil}[1]{\left\lceil #1 \right\rceil}

\catcode`,\active

\catcode`\,12

\theoremstyle{plain}
\newtheorem{thm}{Theorem}[section]

\newtheorem{lem}[thm]{Lemma}
\newtheorem{prop}[thm]{Proposition}

\theoremstyle{definition}

\newtheorem{defn}[thm]{Definition}

\newtheorem{rem}[thm]{Remark}

\theoremstyle{remark}

\setlist[enumerate,1]{leftmargin=2em}

\def\A{\mathbf A}

\def\C{\mathbb C}

\def\H{\mathcal H}
\def\I{\mathbf I}

\def\N{\mathbb N}
\def\P{\mathbf P}
\def\T{\mathbf T}

\def\Z{\mathbb Z}

\def\U{U(\mathfrak{sl}_2)}

\newcommand{\subsetdot}{\subset\mathrel{\mkern-7mu}\mathrel{\cdot}\,}

\title[The Clebsch--Gordan coefficients of $\U$ and Johnson graphs]{The Clebsch--Gordan coefficients of $\U$ and the Terwilliger algebras of Johnson graphs}

\author[Hau-Wen Huang]{Hau-Wen Huang}
\address{
Hau-Wen Huang\\
Department of Mathematics\\
National Central University\\
Chung-Li 32001 Taiwan
}
\email{hauwenh@math.ncu.edu.tw}

\begin{document}
\begin{abstract}
The universal enveloping algebra $\U$ of $\mathfrak{sl}_2$ is a unital associative algebra over $\C$ generated by $E,F,H$ subject to the relations
\begin{align*}
[H,E]=2E,
\qquad 
[H,F]=-2F,
\qquad 
[E,F]=H.
\end{align*}
The element 
$$
\Lambda=EF+FE+\frac{H^2}{2}
$$
is called the Casimir element of $\U$. Let $\Delta:\U\to \U\otimes \U$ denote the comultiplication of $\U$. 
The universal Hahn algebra $\H$ is a unital associative algebra over $\C$ generated by $A,B,C$ and the relations assert that $[A,B]=C$ and each of 
\begin{align*}
[C,A]+2A^2+B,
\qquad 
[B,C]+4BA+2C
\end{align*}
is central in $\H$. 
Inspired by the Clebsch--Gordan coefficients of $\U$, we discover an algebra homomorphism $\natural:\H\to \U\otimes \U$ that maps
\begin{eqnarray*}
A &\mapsto & \frac{H\otimes 1-1\otimes H}{4},
\\
B &\mapsto & \frac{\Delta(\Lambda)}{2},
\\
C &\mapsto & E\otimes F-F\otimes E.
\end{eqnarray*}
By pulling back via $\natural$ any $\U\otimes \U$-module can be considered as an $\H$-module. For any integer $n\geq 0$ there exists a unique $(n+1)$-dimensional irreducible $\U$-module $L_n$ up to isomorphism. We study the decomposition of the $\H$-module $L_m\otimes L_n$ for any integers $m,n\geq 0$. We link these results to the Terwilliger algebras of Johnson graphs. We express the dimensions of the Terwilliger algebras of Johnson graphs in terms of binomial coefficients.
\end{abstract}

\maketitle



\allowdisplaybreaks

\section{Introduction}\label{s:introduction}

Throughout this paper we adopt the following conventions:  
Let $\C$ denote the complex number field. 
Let $\N$ denote the set of nonnegative integers. 
Let $\Z$ denote the ring of integers. 
By an algebra we mean a unital associative algebra. 
By an algebra homomorphism we mean a unital algebra homomorphism. 
By a subalgebra we mean that a nonempty subset is closed under the operations of the parent algebra and has the same unit as the parent algebra. 
In an algebra the bracket $[x,y]$ means the commutator $xy-yx$.  
The unadorned tensor product $\otimes$ is meant to be over $\C$. 
A vacuous summation is interpreted as $0$. A vacuous product is interpreted as $1$. 
For any set $X$ the notation $2^X$ denotes the power set of $X$ and the notation $\subsetdot$ stands for the covering relation of the poset $(2^X,\subseteq)$. 
For any set $X$ the notation $\C^X$ stands for the vector space over $\C$ that has a basis $X$. For any vector space $V$ and any integer $n\geq 1$ the notation $n\cdot V$ stands for 
$
\underbrace{V\oplus V\oplus \cdots \oplus V}_
{\hbox{{\tiny $n$ copies of $V$}}}$.

\begin{defn}\label{defn:U}
The {\it universal enveloping algebra} $\U$ of $\mathfrak{sl}_2$ is an algebra over $\C$ generated by $E,F,H$ subject to the relations
\begin{align}
[H,E]&=2E,
\label{U1}
\\ 
[H,F]&=-2F,
\label{U2}
\\
[E,F]&=H.
\label{U3}
\end{align}
\end{defn}

The algebra $\U$ has a remarkable central element
\begin{gather}\label{Ucasimir}
\Lambda=EF+FE+\frac{H^2}{2}.
\end{gather}
The element $\Lambda$ is called the {\it Casimir element} of $\U$.

\begin{lem}\label{lem:Ln}
For any $n\in \N$ there exists an $(n+1)$-dimensional $\U$-module $L_n$ that has a basis $\{v_i^{(n)}\}_{i=0}^n$ such that 
\begin{align*}
E v_i^{(n)}&=
i v_{i-1}^{(n)}
\quad 
(1\leq i\leq n),
\qquad 
E v_0^{(n)}=0,
\\
F v_i^{(n)}&=
(n-i) v_{i+1}^{(n)}
\quad 
(0\leq i\leq n-1),
\qquad 
F v_n^{(n)}=0,
\\
H v_i^{(n)}&=(n-2i) v_i^{(n)}
\quad 
(0\leq i\leq n).
\end{align*}
\end{lem}
\begin{proof}
It is routine to verify the lemma by using Definition \ref{defn:U}.
\end{proof}

The $\U$-module $L_n$ ($n\in \N$) is the unique $(n+1)$-dimensional irreducible $\U$-module up to isomorphism.

\begin{lem}\label{lem:Casimir&Ln}
For any $n\in \N$ the element $\Lambda$ acts on the $\U$-module $L_n$ as scalar multiplication by $\frac{n(n+2)}{2}$. 
\end{lem}
\begin{proof}
It is straightforward to evaluate the action of $\Lambda$ on the $\U$-module $L_n$ by applying Lemma \ref{lem:Ln} to (\ref{Ucasimir}).
\end{proof}

There exists a unique algebra homomorphism $\Delta:\U\to \U\otimes \U$ given by 
\begin{align}
\Delta(E)&=E\otimes 1+1\otimes E,
\label{DeltaE}
\\
\Delta(F)&=F\otimes 1+1\otimes F,
\label{DeltaF}
\\
\Delta(H)&=H\otimes 1+1\otimes H.
\label{DeltaH}
\end{align}
The map $\Delta$ is called the {\it comultiplication} of $\U$. By pulling back via $\Delta$ any $\U\otimes \U$-module can be considered as a $\U$-module. 
Let $m,n\in \N$ be given. 
By Lemma \ref{lem:Ln} the $\U$-module $L_m\otimes L_n$ has the basis 
\begin{gather}\label{L(mn)basis1}
v_i^{(m)}\otimes v_j^{(n)}
\qquad 
\hbox{($0\leq i\leq m$; $0\leq j\leq n$).}
\end{gather}
The Clebsch--Gordan rule \cite[Proposition V.5.l]{kassel} states that the $\U$-module $L_m\otimes L_n$ is isomorphic to 
$$
\bigoplus\limits_{p=0}^{\min\{m,n\}} L_{m+n-2p}.
$$
Combined with Lemma \ref{lem:Ln} the vectors
\begin{gather}\label{L(mn)basis2}
v_i^{(m+n-2p)}
\qquad 
\hbox{($0\leq p\leq \min\{m,n\}$; $0\leq i\leq m+n-2p$)}
\end{gather}
can be viewed as a basis for $L_m\otimes L_n$. The {\it Clebsch--Gordan coefficients} of $\U$ are the entries of the transition matrix from the basis (\ref{L(mn)basis1}) to the basis (\ref{L(mn)basis2}) for $L_m\otimes L_n$. It is well known that the Clebsch--Gordan coefficients of $\U$ can be expressed in terms of the Hahn polynomials \cite{Klimyk1991,koo1981}. Thus we consider the following algebra arising from the Hahn polynomials:

\begin{defn}
[Definition 1.2, \cite{Huang:even&halved}]
\label{defn:H}
The {\it universal Hahn algebra} $\H$ is an algebra over $\C$ generated by $A,B,C$ and the relations assert that $[A,B]=C$ and each of 
\begin{align}
&[C,A]+2A^2+B,
\label{e:alpha}
\\
&[B,C]+4BA+2C
\label{e:beta}
\end{align}
is central in $\H$.
\end{defn}

Let $\alpha$ and $\beta$ denote the central elements (\ref{e:alpha}) and (\ref{e:beta}) of $\H$ respectively. 
Note that (\ref{e:alpha}) and (\ref{e:beta}) give a central extension of the Askey--Wilson relations corresponding to the Hahn polynomials \cite[Lemma 7.2]{Vidunas:2007}.

By Lemma \ref{lem:Ln} the vectors (\ref{L(mn)basis1}) are the eigenvectors of $H\otimes 1$ and $1\otimes H$ in $L_m\otimes L_n$. By Lemma \ref{lem:Casimir&Ln} the vectors (\ref{L(mn)basis2}) are the eigenvectors of $\Delta(\Lambda)$ in $L_m\otimes L_n$. Based on the above observations we discover the following result:

\begin{thm}\label{thm:H->U2}
There exists a unique algebra homomorphism $\natural:\H\to \U\otimes \U$ that sends
\begin{eqnarray}
 A &\mapsto & \frac{H\otimes 1-1\otimes H}{4}, 
\label{H:A}
\\
 B &\mapsto & \frac{\Delta(\Lambda)}{2}, 
\label{H:B}
\\
C &\mapsto & E\otimes F-F\otimes E.
\label{H:C}
\end{eqnarray}
The algebra homomorphism $\natural$ sends 
\begin{eqnarray}
\alpha 
&\mapsto &
\frac{\Lambda\otimes 1+1\otimes \Lambda}{2}
+
\frac{\Delta(H)^2}{8},
\label{H:alpha}
\\
\beta 
&\mapsto &
\frac{(\Lambda\otimes 1-1\otimes \Lambda)\Delta(H)}{2}.
\label{H:beta}
\end{eqnarray}
\end{thm}

By pulling back via $\natural$ every $\U\otimes \U$-module can be considered as an $\H$-module. 
Let $V$ denote a $\U\otimes \U$-module.  
For any $\theta\in \C$ we define 
$$
V(\theta)=\{v\in V\,|\, \Delta(H)v=\theta v\}.
$$
Recall the element $\Delta(H)\in \U\otimes \U$ from (\ref{DeltaH}).
Clearly $\Delta(H)$ commutes with $H\otimes 1$ and $1\otimes H$.
Since $\Lambda$ is central in $\U$ the element $\Delta(H)$ commutes with 
$\Delta(\Lambda)$. It follows that $\Delta(H)$ commutes with the images of $A$ and $B$ under $\natural$. By Definition \ref{defn:H} the algebra $\H$ is generated by $A$ and $B$. Thus $V(\theta)$ is an $\H$-submodule of $V$ for any $\theta\in \C$.

\begin{lem}
\label{lem:dec_Lmn}
For all $m,n\in \N$ the $\H$-module $L_m\otimes L_n$ is equal to 
\begin{gather*}
\bigoplus_{\ell=0}^{m+n} 
(L_m\otimes L_n)(m+n-2\ell).
\end{gather*}
\end{lem}
\begin{proof}
Immediate from Lemma \ref{lem:Ln}.
\end{proof}

\begin{thm}
\label{thm:Lmn}
Suppose that $m,n\in \N$ and $\ell$ is an integer with $0\leq \ell\leq m+n$. Then the following statements hold:
\begin{enumerate}
\item The $\H$-module $(L_m\otimes L_n)(m+n-2\ell)$ is irreducible with dimension 
$$
\min\{m,\ell\}+\min\{n,\ell\}-\ell+1.
$$

\item Suppose that $m',n'\in \N$ and $\ell'$ is an integer with $0\leq \ell'\leq m'+n'$. Then the $\H$-module $(L_{m'}\otimes L_{n'})(m'+n'-2\ell')$ is isomorphic to $(L_m\otimes L_n)(m+n-2\ell)$ 
if and only if 
\begin{gather}\label{mnl}
(m',n',\ell')\in\{(m,n,\ell),(m+n-\ell,\ell,n),(\ell,m+n-\ell,m),(n,m,m+n-\ell)\}.
\end{gather}
\end{enumerate}
\end{thm}

The outline of this paper is as follows: 
In \S\ref{s:H->U2} we give a proof for Theorem \ref{thm:H->U2}. 
In \S\ref{s:Hmodule} we provide the background for finite-dimensional irreducible $\H$-modules. 
In \S\ref{s:Lmn} we give a proof for Theorem \ref{thm:Lmn}.
In \S\ref{s:Johnson} we connect our results with the Terwilliger algebras of Johnson graphs. In \S\ref{s:dim} we express the dimensions of the Terwilliger algebras of Johnson graphs in terms of binomial coefficients.

\section{Proof for Theorem \ref{thm:H->U2}}\label{s:H->U2}

\begin{lem}\label{lem:H->U2}
The following equations hold in $\U\otimes \U$:
\begin{enumerate}
\item $\Delta(\Lambda)=\Lambda\otimes 1+1\otimes \Lambda+H\otimes H+2(E\otimes F+F\otimes E)$.

\item $
\frac{1}{4}[\Delta(\Lambda),1\otimes H]=F\otimes E-E\otimes F$.

\item $
\frac{1}{4}[\Delta(\Lambda),H\otimes 1]=E\otimes F-F\otimes E$.

\item $
\frac{1}{2}[\Delta(\Lambda),E\otimes F]=(E\otimes F)(1\otimes H-H\otimes 1-2)-H\otimes EF+EF\otimes H$.

\item $\frac{1}{2}[\Delta(\Lambda),F\otimes E]=(F\otimes E)(H\otimes 1-1\otimes H-2)+H\otimes FE-FE\otimes H$.
\end{enumerate}
\end{lem}
\begin{proof}
(i): Applying $\Delta$ to (\ref{Ucasimir}) yields that 
\begin{gather}\label{Delta(Lambda)}
\Delta(\Lambda)=\Delta(EF)+\Delta(FE)+\frac{\Delta(H)^2}{2}.
\end{gather}
Lemma \ref{lem:H->U2}(i) follows by expanding the right-hand side of (\ref{Delta(Lambda)}) by (\ref{DeltaE})--(\ref{DeltaH}) and simplifying the resulting equation by using (\ref{Ucasimir}).

(ii): By Lemma \ref{lem:H->U2}(i) and since $\Lambda$ commutes with $H$ it follows that 
$$
[ \Delta(\Lambda),1\otimes H]=
2([E\otimes F,1\otimes H]+[F\otimes E,1\otimes H]).
$$
Lemma \ref{lem:H->U2}(ii) follows by evaluating the right-hand side of the above equation by using (\ref{U1}) and (\ref{U2}).

(iii): By Lemma \ref{lem:H->U2}(i) and since $\Lambda$ commutes with $H$ it follows that 
$$
[ \Delta(\Lambda),H\otimes 1]=
2([E\otimes F,H\otimes 1]+[F\otimes E,H\otimes 1]).
$$
Lemma \ref{lem:H->U2}(iii) follows by evaluating the right-hand side of the above equation by using (\ref{U1}) and (\ref{U2}).

(iv): By Lemma \ref{lem:H->U2}(i) and since $\Lambda$ commutes with $E$ and $F$ it follows that 
$$
[\Delta(\Lambda),E\otimes F]=
[H\otimes H, E\otimes F]+2[F\otimes E,E\otimes F].
$$
Lemma \ref{lem:H->U2}(iv) follows by evaluating the right-hand side of the above equation by using (\ref{U1})--(\ref{U3}).

(v): By Lemma \ref{lem:H->U2}(i) and since $\Lambda$ commutes with $E$ and $F$ it follows that 
$$
[\Delta(\Lambda),F\otimes E]=
[H\otimes H, F\otimes E]+2[E\otimes F,F\otimes E].
$$
Lemma \ref{lem:H->U2}(v) follows by evaluating the right-hand side of the above equation by using (\ref{U1})--(\ref{U3}). 
\end{proof}

\noindent{\it Proof of Theorem \ref{thm:H->U2}.} 
Let $A^\natural,B^\natural,C^\natural, \alpha^\natural,\beta^\natural$ denote the right-hand sides of (\ref{H:A})--(\ref{H:beta}) respectively. Since $\Lambda$ is central in $\U$ each of $\alpha^\natural$ and $\beta^\natural$ commutes with $A^\natural$ and $B^\natural$. 
To see the existence of $\natural$, by Definition \ref{defn:H} it remains to show that 
\begin{align}
[A^\natural,B^\natural]&=C^\natural,
\label{natural1}
\\
[C^\natural,A^\natural]&=\alpha^\natural-2(A^\natural)^2-B^\natural,
\label{natural2}
\\
[B^\natural,C^\natural] &=\beta^\natural-4B^\natural A^\natural-2 C^\natural.
\label{natural3}
\end{align}

By Lemma \ref{lem:H->U2}(ii), (iii) the equation (\ref{natural1}) follows. 
Using (\ref{U1}) and (\ref{U2}) yields that 
\begin{align*}
[E\otimes F-F\otimes E, H\otimes 1]&=-2(E\otimes F+F\otimes E),
\\
[E\otimes F-F\otimes E, 1\otimes H]&=2(E\otimes F+F\otimes E).
\end{align*} 
Hence the left-hand side of (\ref{natural2}) is equal to 
$
-(E\otimes F+F\otimes E)$. 
Using Lemma \ref{lem:H->U2}(i) it is routine to verify that the right-hand side of (\ref{natural2}) is also equal to $
-(E\otimes F+F\otimes E)$. The equation (\ref{natural2}) follows. By Lemma \ref{lem:H->U2}(iv), (v) the left-hand side of (\ref{natural3}) is equal to $(E\otimes F+F\otimes E)(1\otimes H-H\otimes 1)+2(F\otimes E-E\otimes F)$ plus 
\begin{gather}\label{natural3'}
(EF+FE)\otimes H-H\otimes (EF+FE).
\end{gather}
By Lemma \ref{lem:H->U2}(i) the right-hand side of (\ref{natural3}) is equal to $(E\otimes F+F\otimes E)(1\otimes H-H\otimes 1)+2(F\otimes E-E\otimes F)$ plus 
\begin{gather}\label{natural3''}
\frac{H\otimes H}{2}
(1\otimes H-H\otimes 1)
+\Lambda\otimes H
-H\otimes \Lambda.
\end{gather}
Substituting (\ref{Ucasimir}) into (\ref{natural3''}) yields that (\ref{natural3'}) is equal to (\ref{natural3''}). The equation (\ref{natural3}) follows. Since (\ref{natural1})--(\ref{natural3}) hold the existence of $\natural$ follows. 
Since the algebra $\H$ is generated by $A,B,C$ the uniqueness of $\natural$ follows. 
\hfill $\square$

\section{Preliminaries on the finite-dimensional irreducible $\H$-modules}\label{s:Hmodule}

In this section we state the preliminaries on the finite-dimensional irreducible $\H$-modules. Here we omit the proofs because they are analogous to the proofs in the literature \cite{Huang:2015,SH:2019-1,Huang:BImodule}. 

\begin{prop}\label{prop:Vd(a,b)}
For any $a,b\in \C$ and $d\in \N$ there exists a $(d+1)$-dimensional $\H$-module $V_d(a,b)$ satisfying the following conditions:
\begin{enumerate}
\item There exists a basis $\{v_i\}_{i=0}^d$ for $V_d(a,b)$ such that 
\begin{align*}
(A-\theta_i) v_i
&=
\varphi_i v_{i-1}
\quad 
(1\leq i\leq d),
\qquad 
(A-\theta_0)v_0=0,
\\
(B-\theta_i^*) v_i
&=
v_{i+1}
\quad 
(0\leq i\leq d-1),
\qquad 
(B-\theta_d^*) v_d
=0,
\end{align*}
where 
\begin{align*}
\theta_i&=
\frac{a+d}{2}-i
\qquad 
(0\leq i\leq d),
\\
\theta_i^* &=(b+i)(b+i+1)
\qquad 
(0\leq i\leq d),
\\
\varphi_i &=
i(i-d-1)
(a-b-i)
\qquad 
(1\leq i\leq d).
\end{align*}

\item The elements $\alpha$ and $\beta$ act on $V_d(a,b)$ as scalar multiplication by \begin{align*}
\eta &=
\frac{a^2+d(d+2)}{2}+b(b+d+1),
\\
\eta^* &=
2ab(b+d+1),
\end{align*} 
respectively.
\end{enumerate}
\end{prop}

For any $a,b\in \C$ and $d\in \N$ the notation $V_d(a,b)$ will always stand for the $\H$-module given in Proposition \ref{prop:Vd(a,b)}.

\begin{prop}\label{prop:Vd_universal}
Let the notations be as in Proposition \ref{prop:Vd(a,b)}.
Suppose that $V$ is an $\H$-module satisfying the following conditions:
\begin{enumerate}
\item There is a vector $v\in V$ satisfying
\begin{gather*}
Av=\theta_0 v,
\\
(A-\theta_1)(B-\theta_0^*)v=\varphi_1 v,
\\
\alpha v=\eta v, 
\qquad 
\beta v=\eta^* v.
\end{gather*}

\item $\frac{a-d}{2}-1$ is not an eigenvalue of $A$ in $V$.
\end{enumerate}
Then there exists a unique $\H$-module homomorphism 
$
V_d(a, b) \to V
$ that sends $v_0$ to $v$.
\end{prop}

\begin{thm}\label{thm:Vd_irr}
For any $a,b\in \C$ and $d\in \N$ the $\H$-module $V_d(a,b)$ is irreducible if and only if 
$$
a-b,-a-b\not\in\{1,2,\ldots,d\}.
$$
\end{thm}

By an argument similar to the proofs of \cite[Theorem 4.7]{Huang:2015}, \cite[Theorem 6.3]{Huang:BImodule} and \cite[Theorem 6.3]{SH:2019-1}, it can be shown that for any finite-dimensional irreducible $\H$-module $V$ there exist $a,b\in \C$ and $d\in \N$ such that the $\H$-module $V$ is isomorphic to $V_d(a,b)$. Moreover the parameters $a,b,d$ can be determined in the following way:

Recall from Schur's lemma that every central element of an algebra $\mathcal A$ over $\C$ acts on a finite-dimensional irreducible $\mathcal A$-module as scalar multiplication. 

\begin{thm}\label{thm:Vd_onto}
Suppose that $V$ is a finite-dimensional irreducible $\H$-module.  
Let ${\rm tr}\, A$ denote the trace of $A$ on $V$. 
Let $\eta$ denote the unique eigenvalue of $\alpha$ on $V$.
For any $a,b\in \C$ and $d\in \N$, the $\H$-module $V$ is isomorphic to $V_d(a,b)$
if and only if the following conditions hold:
\begin{enumerate}
\item $d+1$ is equal to the dimension of $V$.

\item $a(d+1)=2\cdot {\rm tr}\, A$. 

\item $b$ is a root of 
$x^2+(d+1) x+\frac{a^2+d(d+2)}{2}-\eta$. 
\end{enumerate} 
\end{thm}

\section{Proof for Theorem \ref{thm:Lmn}}\label{s:Lmn}

Recall the $\U$-module $L_n$ ($n\in \N$) from Lemma \ref{lem:Ln}.

\begin{lem}\label{lem:EFH_Lmn}
Suppose that $m,n\in \N$. For any integers $i,j$ with $0\leq i\leq m$ and $0\leq j\leq n$ the following equations hold on the $\U\otimes \U$-module $L_m\otimes L_n$:
\begin{enumerate}
\item 
$
\Delta(E)\, v_i^{(m)}\otimes v_j^{(n)}
=
i v_{i-1}^{(m)}\otimes v_j^{(n)}
+
j v_i^{(m)}\otimes v_{j-1}^{(n)}.
$

\item 
$
\Delta(F)\, v_i^{(m)}\otimes v_j^{(n)}
=
(m-i) v_{i+1}^{(m)}\otimes v_j^{(n)}
+
(n-j) v_i^{(m)}\otimes v_{j+1}^{(n)}.
$
\end{enumerate}
Here $v_{-1}^{(m)}, v_{m+1}^{(m)}$ are interpreted as any vectors of $L_m$ and $v_{-1}^{(n)}, v_{n+1}^{(n)}$ are interpreted as any vectors of $L_n$.
\end{lem}
\begin{proof}
It is straightforward to verify the lemma by using (\ref{DeltaE}), (\ref{DeltaF}) and Lemma \ref{lem:Ln}.
\end{proof}

Recall the notation $V(\theta)$ for any $\U\otimes \U$-module $V$ and any scalar $\theta\in \C$ from below Theorem \ref{thm:H->U2}.

\begin{thm}\label{thm:LHmodule}
Suppose that $m,n\in \N$ and $\ell$ is an integer with $0\leq \ell\leq m+n$. Then the $\H$-module $(L_m\otimes L_n)(m+n-2\ell)$ is isomorphic to the irreducible $\H$-module $V_d(a,b)$ where 
\begin{align}
a&=\min\{n,\ell\}-\min\{m,\ell\}+\frac{m-n}{2},
\label{LH:a}
\\
b&=-\frac{m+n}{2}-1,
\label{LH:b}
\\
d &=\min\{m,\ell\}+\min\{n,\ell\}-\ell.
\label{LH:d}
\end{align} 
\end{thm}
\begin{proof}
Let $a,b,d$ denote the current values (\ref{LH:a})--(\ref{LH:d}) respectively. 
Let $\{v_i\}_{i=0}^d$ denote the basis for $V_d(a,b)$ given in Proposition \ref{prop:Vd(a,b)}(i). 
Let $\{\theta_i\}_{i=0}^d,\{\theta_i^*\}_{i=0}^d,\{\varphi_i\}_{i=1}^d,\eta,\eta^*$ denote the parameters described in Proposition \ref{prop:Vd(a,b)}. 
Applying Theorem \ref{thm:Vd_irr} it is routine to verify the irreducibility of the $\H$-module $V_d(a,b)$ by dividing the argument into the following cases: 
\begin{enumerate}
\item[(a)] $\ell<\min\{m,n\}$; 
\item[(b)] $m\leq \ell \leq n$; 
\item[(c)] $n\leq \ell\leq m$; 
\item[(d)] $\ell>\max\{m,n\}$. 
\end{enumerate}
To see Theorem \ref{thm:LHmodule} it remains to show that the $\H$-module $(L_m\otimes L_n)(m+n-2\ell)$ is isomorphic to $V_d(a,b)$.

Recall the element $\Delta(H)\in \U\otimes \U$ from (\ref{DeltaH}). By Lemma \ref{lem:Ln} the vectors
\begin{gather}\label{uk}
v_{\ell-k}^{(m)}  \otimes v_{k}^{(n)}
\qquad (\max\{m,\ell\}-m\leq k\leq \min\{n,\ell\})
\end{gather}
are a basis for $(L_m\otimes L_n)(m+n-2\ell)$. Replacing the index $k$ by $\min\{n,\ell\}-i$ the vectors (\ref{uk}) can be written as 
\begin{gather*}
u_i=v_{\ell-\min\{n,\ell\}+i}^{(m)}  \otimes v_{\min\{n,\ell\}-i}^{(n)}
\qquad (0\leq i\leq d).
\end{gather*}
Hence $(L_m\otimes L_n)(m+n-2\ell)$ has dimension $d+1$.
Recall the image of $A$ under $\natural$ from (\ref{H:A}). 
Using Lemma \ref{lem:Ln} it is straightforward to check that 
\begin{gather}\label{Aui}
A u_i =\theta_i u_i
\qquad 
(0\leq i\leq d).
\end{gather}
Hence $\{\theta_i\}_{i=0}^d$ are all eigenvalues of $A$ in $(L_m\otimes L_n)(m+n-2\ell)$. 
Using Lemma \ref{lem:EFH_Lmn} yields that 
\begin{align*}
\Delta(EF) u_0 &=
c_0 u_0
+
c_1 u_1,
\\
\Delta(FE) u_0 &=
c_2 u_0
+
c_1 u_1,
\end{align*}
where 
\begin{align*}
c_0 &=
(\min\{n,\ell\}+1)(n-\min\{n,\ell\})
+
(\ell-\min\{n,\ell\}+1)(m-\ell+\min\{n,\ell\}),
\\
c_1 &=\min\{n,\ell\}(m-\ell+\min\{n,\ell\}),
\\
c_2 &=\min\{n,\ell\}(n-\min\{n,\ell\}+1)
+
(\ell-\min\{n,\ell\})(m-\ell+\min\{n,\ell\}+1),
\end{align*}
and $u_1$ is interpreted as any vector of $(L_m\otimes L_n)(m+n-2\ell)$ when $d=0$.
In each of the cases (a)--(d) the following equalities hold:
\begin{align*}
\varphi_1&=-c_1,
\\
\theta_0^*+\varphi_1&=\frac{c_0+c_2}{2}+\frac{(m+n-2\ell)^2}{4}.
\end{align*} 
Here $\varphi_1$ is interpreted as zero when $d=0$. 
Recall the image of $B$ under $\natural$ from (\ref{H:B}). 
Combined with (\ref{Ucasimir}) this yields that 
$$(B-\theta_0^*) u_0=
\varphi_1 
(u_0-u_1).
$$
Applying (\ref{Aui}) to the above equality yields that 
$$
(A-\theta_1)(B-\theta_0^*) 
u_0
=
\varphi_1
u_0.
$$ 
Recall the images of $\alpha$ and $\beta$ under $\natural$ from (\ref{H:alpha}) and (\ref{H:beta}). 
By Lemma \ref{lem:Casimir&Ln} the elements $\alpha$ and $\beta$ act on $(L_m\otimes L_n)(m+n-2\ell)$ as scalar multiplication by 
\begin{gather}
\frac{m(m+2)+n(n+2)}{4}+\frac{(m+n-2\ell)^2}{8},
\label{eta}
\\
\frac{m(m+2)-n(n+2)}{4}(m+n-2\ell),
\label{etas}
\end{gather}
respectively. In each of the cases (a)--(d)  the scalars $\eta$ and $\eta^*$ are equal to 
(\ref{eta}) and (\ref{etas}) respectively.

In light of the above results, it follows from Proposition \ref{prop:Vd_universal} that there exists a unique $\H$-module homomorphism 
\begin{gather}\label{V->LH}
V_d(a,b)\to (L_m\otimes L_n)(m+n-2\ell)
\end{gather} 
that sends $v_0$ to $u_0$.
Since the $\H$-module $V_d(a,b)$ is irreducible, it follows that (\ref{V->LH}) is injective. 
Since the dimensions of $(L_m\otimes L_n)(m+n-2\ell)$ and $V_d(a,b)$ are equal to $d+1$, it follows that (\ref{V->LH}) is an isomorphism.  The result follows.
\end{proof}

\noindent
{\it Proof of Theorem \ref{thm:Lmn}}. 
(i): Immediate from Theorem \ref{thm:LHmodule}.

(ii): Let $a,b,d$ denote the parameters (\ref{LH:a})--(\ref{LH:d}) respectively. Let $a',b',d'$ denote the right-hand sides of (\ref{LH:a})--(\ref{LH:d}) with the replacement of $(m,n,\ell)$ by $(m',n',\ell')$ respectively. 
By Theorem \ref{thm:LHmodule} the $\H$-modules $(L_m\otimes L_n)(m+n-2\ell)$ and $(L_{m'}\otimes L_{n'})(m'+n'-2\ell')$ are isomorphic to the finite-dimensional irreducible $\H$-modules $V_d(a,b)$ and $V_{d'}(a',b')$ respectively. Note that $b\leq -(d+1)$ and $b'\leq -(d'+1)$. Combined with Theorem \ref{thm:Vd_onto} the $\H$-modules $V_d(a,b)$ and $V_{d'}(a',b')$ are isomorphic if and only if 
\begin{gather}\label{abd}
(a,b,d)=(a',b',d').
\end{gather} 
To see Theorem \ref{thm:Lmn}(ii) it suffices to show that (\ref{abd}) is equivalent to the condition (\ref{mnl}). By symmetry we may divide the argument into the following ten cases:

(a): Suppose that $\ell< \min\{m,n\}$ and $\ell'< \min\{m',n'\}$. The condition (\ref{mnl}) turns out $(m',n',\ell')=(m,n,\ell)$. In this case 
\begin{align*}
(a,b,d)&=\left(
\frac{m-n}{2},-\frac{m+n}{2}-1,\ell
\right),
\\
(a',b',d')&=\left(
\frac{m'-n'}{2},-\frac{m'+n'}{2}-1,\ell'
\right).
\end{align*}
Solving (\ref{abd}) yields that $(m',n',\ell')=(m,n,\ell)$.

(b): Suppose that $\ell<\min\{m,n\}$ and $m'\leq \ell'\leq n'$. The condition (\ref{mnl}) turns out $(m',n',\ell')=(\ell,m+n-\ell,m)$. In this case 
\begin{align*}
(a,b,d)&=\left(
\frac{m-n}{2},-\frac{m+n}{2}-1,\ell
\right),
\\
(a',b',d')&=\left(
\ell'-\frac{m'+n'}{2},-\frac{m'+n'}{2}-1,m'
\right).
\end{align*}
Solving (\ref{abd}) yields that $(m',n',\ell')=(\ell,m+n-\ell,m)$.

(c): Suppose that $\ell<\min\{m,n\}$ and $n'\leq \ell'\leq m'$. The condition (\ref{mnl}) turns out $(m',n',\ell')=(m+n-\ell,\ell,n)$. In this case 
\begin{align*}
(a,b,d)&=\left(
\frac{m-n}{2},-\frac{m+n}{2}-1,\ell
\right),
\\
(a',b',d')&=\left(
\frac{m'+n'}{2}-\ell',-\frac{m'+n'}{2}-1,n'
\right).
\end{align*}
Solving (\ref{abd}) yields that $(m',n',\ell')=(m+n-\ell,\ell,n)$.

(d): Suppose that $\ell<\min\{m,n\}$ and $\ell'>\max\{m',n'\}$. The condition (\ref{mnl}) turns out $(m',n',\ell')=(n,m,m+n-\ell)$. In this case 
\begin{align*}
(a,b,d)&=\left(
\frac{m-n}{2},-\frac{m+n}{2}-1,\ell
\right),
\\
(a',b',d')&=
\left(
\frac{n'-m'}{2},-\frac{m'+n'}{2}-1,m'+n'-\ell'
\right).
\end{align*}
Solving (\ref{abd}) yields that $(m',n',\ell')=(n,m,m+n-\ell)$.

(e): Suppose that $m\leq \ell\leq n$ and $m'\leq \ell'\leq n'$. The condition (\ref{mnl}) turns out $(m',n',\ell')=(m,n,\ell)$. In this case 
\begin{align*}
(a,b,d)&=\left(
\ell-\frac{m+n}{2},-\frac{m+n}{2}-1,m
\right),
\\
(a',b',d')&=
\left(
\ell'-\frac{m'+n'}{2},-\frac{m'+n'}{2}-1,m'
\right).
\end{align*}
Solving (\ref{abd}) yields that $(m',n',\ell')=(m,n,\ell)$.

(f): Suppose that $m\leq \ell\leq n$ and $n'\leq \ell'\leq m'$. The condition (\ref{mnl}) turns out $(m',n',\ell')=(n,m,m+n-\ell)$. In this case 
\begin{align*}
(a,b,d)&=\left(
\ell-\frac{m+n}{2},-\frac{m+n}{2}-1,m
\right),
\\
(a',b',d')&=
\left(
\frac{m'+n'}{2}-\ell',-\frac{m'+n'}{2}-1,n'
\right).
\end{align*}
Solving (\ref{abd}) yields that $(m',n',\ell')=(n,m,m+n-\ell)$.

(g): Suppose that $m\leq \ell\leq n$ and $\ell'>\max\{m',n'\}$. The condition (\ref{mnl}) turns out $(m',n',\ell')=(m+n-\ell,\ell,n)$. In this case 
\begin{align*}
(a,b,d)&=\left(
\ell-\frac{m+n}{2},-\frac{m+n}{2}-1,m
\right),
\\
(a',b',d')&=
\left(
\frac{n'-m'}{2},-\frac{m'+n'}{2}-1,m'+n'-\ell'
\right).
\end{align*}
Solving (\ref{abd}) yields that $(m',n',\ell')=(m+n-\ell,\ell,n)$.

(h): Suppose that $n\leq \ell\leq m$ and $n'\leq \ell'\leq m'$. The condition (\ref{mnl}) turns out $(m',n',\ell')=(m,n,\ell)$. In this case 
\begin{align*}
(a,b,d)&=\left(
\frac{m+n}{2}-\ell,-\frac{m+n}{2}-1,n
\right),
\\
(a',b',d')&=
\left(
\frac{m'+n'}{2}-\ell',-\frac{m'+n'}{2}-1,n'
\right).
\end{align*}
Solving (\ref{abd}) yields that $(m',n',\ell')=(m,n,\ell)$.

(i): Suppose that $n\leq \ell\leq m$ and $\ell'>\max\{m',n'\}$. The condition (\ref{mnl}) turns out $(m',n',\ell')=(\ell,m+n-\ell,m)$. In this case 
\begin{align*}
(a,b,d)&=\left(
\frac{m+n}{2}-\ell,-\frac{m+n}{2}-1,n
\right),
\\
(a',b',d')&=
\left(
\frac{n'-m'}{2},-\frac{m'+n'}{2}-1,m'+n'-\ell'
\right).
\end{align*}
Solving (\ref{abd}) yields that $(m',n',\ell')=(\ell,m+n-\ell,m)$.

(j): Suppose that $\ell> \max\{m,n\}$ and $\ell'>\max\{m',n'\}$. The condition (\ref{mnl}) turns out $(m',n',\ell')=(m,n,\ell)$. In this case 
\begin{align*}
(a,b,d)&=\left(
\frac{n-m}{2},-\frac{m+n}{2}-1,m+n-\ell
\right),
\\
(a',b',d')&=
\left(
\frac{n'-m'}{2},-\frac{m'+n'}{2}-1,m'+n'-\ell'
\right).
\end{align*}
Solving (\ref{abd}) yields that $(m',n',\ell')=(m,n,\ell)$. Based on the above results (a)--(j) the statement (ii) of Theorem \ref{thm:Lmn} follows.
\hfill $\square$

\section{The connection to the Terwilliger algebras of Johnson graphs}\label{s:Johnson}

In the rest of this paper the notation $\Omega$ stands for a finite set with size $D$. 
In the study of the Terwilliger algebra of the hypercube \cite{hypercube2002}, it was mentioned the following connection of $\U$ and $2^\Omega$:

\begin{thm}
[Theorem 13.2, \cite{hypercube2002}]
\label{thm:U&2Omega}
There exists a $\U$-module structure on $\C^{2^\Omega}$ given by 
\begin{align*}
E x &=\sum_{y\subsetdot x} y
\qquad 
\hbox{for all $x\in 2^\Omega$},
\\ 
F x &=\sum_{x\subsetdot y} y
\qquad 
\hbox{for all $x\in 2^\Omega$},
\\ 
H x &= (D-2|x|) x
\qquad 
\hbox{for all $x\in 2^\Omega$},
\end{align*} 
where the first summation is over all $y\in 2^\Omega$ with $y\subsetdot x$ and the second summation is over all $y\in 2^\Omega$ with $x\subsetdot y$.
\end{thm}

\begin{lem}
\label{lem:Lambda_2Omega}
The action of $\Lambda$ on the $\U$-module $\C^{2^\Omega}$ is as follows:
\begin{align*}
\Lambda x=
\left(
D+
\frac{(D-2|x|)^2}{2}
\right)
x
+
2\sum_{\substack{|y|=|x|\\ x\cap y\subsetdot x}} y
\qquad 
\hbox{for all $x\in 2^\Omega$},
\end{align*}
where the summation is over all $y\in 2^\Omega$ with $|y|=|x|$ and $x\cap y\subsetdot x$.
\end{lem}
\begin{proof}
Using Theorem \ref{thm:U&2Omega} yields that 
\begin{align*}
EF x &=
(D-|x|)x+\sum_{\substack{|y|=|x|\\ x\cap y\subsetdot x}} y
\qquad 
\hbox{for all $x\in 2^\Omega$},
\\
FE x &=
|x|x+\sum_{\substack{|y|=|x|\\ x\cap y\subsetdot x}} y
\qquad 
\hbox{for all $x\in 2^\Omega$},
\\ 
H^2 x &= (D-2|x|)^2 x 
\qquad 
\hbox{for all $x\in 2^\Omega$}.
\end{align*}
Combined with (\ref{Ucasimir}) this yields the lemma.
\end{proof}

Recall that the binomial coefficients are defined by
\begin{gather*}
{n\choose \ell}=\prod_{i=1}^\ell \frac{n-i+1}{\ell-i+1}
\qquad 
\hbox{for all $\ell,n\in \N$}.
\end{gather*}
For notational convenience let
$$
m_i(n)=\frac{n-2i+1}{n-i+1} 
{n\choose i}
$$
for all integers $i,n$ with $0\leq i\leq \floor{\frac{n}{2}}$.

\begin{thm}
[Theorem 10.2, \cite{hypercube2002}]
\label{thm:dec2Omega}
The $\U$-module $\C^{2^\Omega}$ is isomorphic to 
\begin{gather*}
\bigoplus_{i=0}^{\floor{\frac{D}{2}}}
m_i(D)
\cdot 
L_{D-2i}.
\end{gather*}
\end{thm}

Fix an element $x_0\in 2^\Omega$. By Theorem \ref{thm:U&2Omega}
the spaces $\C^{2^{\Omega\setminus x_0}}$ and $\C^{2^{x_0}}$ are $\U$-modules. Hence $\C^{2^{\Omega\setminus x_0}}\otimes \C^{2^{x_0}}$ has a $\U\otimes \U$-module structure.  We consider the linear transformation 
$\iota(x_0):\C^{2^\Omega}\to \C^{2^{\Omega\setminus x_0}}\otimes \C^{2^{x_0}}$ given by 
\begin{eqnarray*}
x &\mapsto & (x\setminus x_0)\otimes (x\cap x_0)
\qquad 
\hbox{for all $x\in 2^\Omega$}.
\end{eqnarray*}
Note that $\iota(x_0)$ is an isomorphism.

\begin{lem}\label{lem:X&DeltaX}
For any $x_0\in 2^\Omega$ and $X\in \U$ the following diagram commutes:
\begin{table}[H]
\centering
\begin{tikzpicture}
\matrix(m)[matrix of math nodes,
row sep=4em, column sep=4em,
text height=1.5ex, text depth=0.25ex]
{
\C^{2^\Omega}
&\C^{2^{\Omega\setminus x_0}}\otimes \C^{2^{x_0}}\\
\C^{2^\Omega}
&\C^{2^{\Omega\setminus x_0}}\otimes \C^{2^{x_0}}\\
};
\path[->,font=\scriptsize,>=angle 90]
(m-1-1) edge node[left] {$X$} (m-2-1)
(m-1-1) edge node[above] {$\iota(x_0)$} (m-1-2)
(m-2-1) edge node[below] {$\iota(x_0)$} (m-2-2)
(m-1-2) edge node[right] {$\Delta(X)$} (m-2-2);
\end{tikzpicture}
\end{table}
\end{lem}
\begin{proof}
Let $x_0\in 2^\Omega$ be given. 
Using (\ref{DeltaE})--(\ref{DeltaH}) and Theorem \ref{thm:U&2Omega} it is routine to verify that the  diagram commutes for $X=E,F,H$. Since the algebra $\U$ is generated by $E,F,H$ the lemma follows.
\end{proof}

By identifying $\C^{2^\Omega}$ with $\C^{2^{\Omega\setminus x_0}}\otimes \C^{2^{x_0}}$ via $\iota(x_0)$, this induces a $\U\otimes \U$-module structure on $\C^{2^\Omega}$ and we denote by $\C^{2^\Omega}(x_0)$.

\begin{lem}
\label{lem:U2&2Omega}
For any $x_0\in 2^\Omega$ the $\U\otimes \U$-module $\C^{2^\Omega}(x_0)$ is isomorphic to 
\begin{gather*}
\bigoplus_{i=0}^{\floor{\frac{D-|x_0|}{2}}}
\bigoplus_{j=0}^{\floor{\frac{|x_0|}{2}}}
m_i(D-|x_0|)
\,
m_j(|x_0|)
\cdot 
L_{D-|x_0|-2i}\otimes L_{|x_0|-2j}.
\end{gather*}
\end{lem}
\begin{proof}
Immediate from Theorem \ref{thm:dec2Omega} and the distributive law of $\otimes$ over $\oplus$.
\end{proof}

By pulling back via $\natural$ the $\U\otimes \U$-module $\C^{2^\Omega}(x_0)$ is an $\H$-module. Suppose that $k$ is an integer with $0\leq k\leq D$. We adopt the notation ${\Omega \choose k}$ to denote the set of all $k$-element subsets of $\Omega$. Applying Lemma \ref{lem:X&DeltaX} with $X=H$ the action of $\Delta(H)$ on the $\U\otimes \U$-module $\C^{2^\Omega}(x_0)$ is identical to the action of $H$ on the $\U$-module on $\C^{2^\Omega}$. Combined with Theorem \ref{thm:U&2Omega} this yields that 
\begin{gather}\label{C(Omega k)}
\C^{\Omega\choose k}=\C^{2^\Omega}(x_0)(D-2k).
\end{gather}
It follows that $\C^{{\Omega \choose k}}$ is an $\H$-submodule of $\C^{2^\Omega}(x_0)$ and we denote this $\H$-module by $\C^{{\Omega \choose k}}(x_0)$.

For notational convenience let 
\begin{gather}\label{P}
\P(k)=
\left\{
(i,j)\in \Z^2
\,\Big|\,
0\leq i\leq \frac{D-k}{2},
0\leq j\leq \min
\left\{
D-k-i,k-i,\frac{k}{2}
\right\} 
\right\}.
\end{gather}

\begin{lem}
\label{lem:H&Omega k}
Suppose that $k$ is an integer with $0\leq k\leq D$. 
For any $x_0\in {\Omega\choose k}$ the $\H$-module $\C^{{\Omega \choose k}}(x_0)$ is isomorphic to 
\begin{align*}
\bigoplus_{(i,j)\in \P(k)}
m_i(D-k)\,
m_j(k)
\cdot 
(L_{D-k-2i}\otimes L_{k-2j})(D-2k).
\end{align*}
Moreover the $\H$-module $(L_{D-k-2i}\otimes L_{k-2j})(D-2k)$ is irreducible for each $(i,j)\in \P(k)$.
\end{lem}
\begin{proof}
In view of Lemma \ref{lem:U2&2Omega} and (\ref{C(Omega k)}) the $\H$-module $\C^{\Omega\choose k}(x_0)$ is isomorphic to 
$$
\bigoplus_{i=0}^{\floor{\frac{D-k}{2}}}
\bigoplus_{j=0}^{\floor{\frac{k}{2}}}
m_i(D-k)
\,
m_j(k)
\cdot 
(L_{D-k-2i}\otimes L_{k-2j})(D-2k).
$$
Applying Lemma \ref{lem:dec_Lmn} with $(m,n,\ell)=(D-k-2i,k-2j,k-i-j)$ yields that $(L_{D-k-2i}\otimes L_{k-2j})(D-2k)=\{0\}$ for any integer $j\not\in \left\{0,1,\ldots,\min\{D-k-i,k-i,\floor{\frac{k}{2}}\}\right\}$. Combined with Theorem \ref{thm:Lmn}(i) the lemma follows.
\end{proof}

\begin{thm}
\label{thm:H&Omega k}
Suppose that $k$ is an integer with $0\leq k\leq D$. 
For any $x_0\in {\Omega\choose k}$ the following statements hold:
\begin{enumerate}
\item Suppose that $k\not=\frac{D}{2}$. Then the $\H$-module $\C^{{\Omega \choose k}}(x_0)$ is isomorphic to 
\begin{align*}
\bigoplus_{(i,j)\in \P(k)}
m_i(D-k)
m_j(k)
\cdot 
(L_{D-k-2i}\otimes L_{k-2j})(D-2k).
\end{align*}
Moreover the irreducible $\H$-modules $(L_{D-k-2i}\otimes L_{k-2j})(D-2k)$ for all $(i,j)\in \P(k)$ are mutually non-isomorphic.

\item Suppose that $k=\frac{D}{2}$. Then the $\H$-module $\C^{{\Omega \choose k}}(x_0)$ is isomorphic to a direct sum of 
\begin{align*}
\bigoplus_{i=0}^{\floor{\frac{D}{4}}}
m_i\left(\frac{D}{2}\right)^2 
\cdot 
(L_{\frac{D}{2}-2i}\otimes L_{\frac{D}{2}-2i})(0)
\end{align*}
and 
\begin{align*}
\bigoplus_{i=0}^{\floor{\frac{D}{4}}}
\bigoplus_{j=i+1}^{\floor{\frac{D}{4}}}
2
m_i\left(\frac{D}{2}\right)
m_j\left(\frac{D}{2}\right)
\cdot 
(L_{\frac{D}{2}-2i}\otimes L_{\frac{D}{2}-2j})(0).
\end{align*}
Moreover the irreducible $\H$-modules $(L_{\frac{D}{2}-2i}\otimes L_{\frac{D}{2}-2i})(0)$ for all integers $i$ with $0\leq i\leq \floor{\frac{D}{4}}$ and the irreducible $\H$-modules $(L_{\frac{D}{2}-2i}\otimes L_{\frac{D}{2}-2j})(0)$ for all integers $i,j$ with $0\leq i\leq \floor{\frac{D}{4}}$ and $i+1\leq j\leq \floor{\frac{D}{4}}$ are mutually non-isomorphic.
\end{enumerate}
\end{thm}
\begin{proof}
(i): The first assertion is immediate from Lemma \ref{lem:H&Omega k}. Suppose that there are $(i,j),(i',j')\in \P(k)$ such that the $\H$-modules $(L_{D-k-2i}\otimes L_{k-2j})(D-2k)$ and $(L_{D-k-2i'}\otimes L_{k-2j'})(D-2k)$ are isomorphic. Applying Theorem \ref{thm:Lmn}(ii) with $(m,n,\ell)=(D-k-2i,k-2j,k-i-j)$ and $(m',n',\ell')=(D-k-2i',k-2j',k-i'-j')$ yields that $(i,j)=(i',j')$. The second assertion follows.

(ii): In this case
$
\P(k)=
\textstyle
\left\{
(i,j)\in \Z^2\,|\,0\leq i,j\leq \floor{\frac{D}{4}}
\right\}$.
It follows from Lemma \ref{lem:H&Omega k} that the $\H$-module $\C^{{\Omega\choose k}}(x_0)$ is isomorphic to 
\begin{gather*}
\bigoplus_{i=0}^{\floor{\frac{D}{4}}}
\bigoplus_{j=0}^{\floor{\frac{D}{4}}}
m_i\left(\frac{D}{2}\right)
m_j\left(\frac{D}{2}\right)
\cdot 
(
L_{\frac{D}{2}-2i}\otimes L_{\frac{D}{2}-2j}
)(0).
\end{gather*}
Let $(i,j),(i',j')\in \P(k)$ be given. Applying Theorem \ref{thm:Lmn}(ii) with $(m,n,\ell)=(\frac{D}{2}-2i,\frac{D}{2}-2j,\frac{D}{2}-i-j)$ and $(m',n',\ell')=(\frac{D}{2}-2i',\frac{D}{2}-2j',\frac{D}{2}-i'-j')$ yields that the $\H$-modules $(L_{\frac{D}{2}-2i}\otimes L_{\frac{D}{2}-2j})(0)$ and $(L_{\frac{D}{2}-2i'}\otimes L_{\frac{D}{2}-2j'})(0)$ are isomorphic if and only if $(i',j')\in\{(i,j), (j,i)\}$. 
By the above comments the statement (ii) follows.
\end{proof}

Recall that the algebra $\H$ is generated by $A$ and $B$. A description for the $\H$-module $\C^{2^\Omega}(x_0)$ is as follows:

\begin{thm}\label{thm:H&2Omega}
For any $x_0\in 2^\Omega$ the actions of $A$ and $B$ on the $\H$-module $\C^{2^\Omega}(x_0)$ are as follows:
\begin{align}
A x &=
\left(
\frac{D}{4}-\frac{|x_0\setminus x|+|x\setminus x_0|}{2}
\right) 
x
\qquad 
\hbox{for all $x\in 2^\Omega$},
\label{A_2Omega(x0)}
\\
B x &=
\left(
\frac{D}{2}
+
\frac{(D-2|x|)^2}{4}
\right)
x
+
\sum_{\substack{|y|=|x|\\ x\cap y\subsetdot x}} y
\qquad 
\hbox{for all $x\in 2^\Omega$}.
\label{B_2Omega(x0)}
\end{align}
\end{thm}
\begin{proof}
Let $x_0\in 2^\Omega$ be given. Recall the linear isomorphism $\iota(x_0):\C^{2^\Omega}\to \C^{2^{\Omega\setminus x_0}}\otimes \C^{2^{x_0}}$ from above Lemma \ref{lem:X&DeltaX}.
Using Theorem \ref{thm:U&2Omega} yields that 
\begin{align*}
(H\otimes 1)x &=
(D-|x_0|-2|x\setminus x_0|) x 
\qquad 
\hbox{for all $x\in 2^\Omega$},
\\
(1\otimes H)x &=
(|x_0|-2|x\cap x_0|) x 
\qquad 
\hbox{for all $x\in 2^\Omega$}.
\end{align*}
Combined with (\ref{H:A}) this yields (\ref{A_2Omega(x0)}). 
Applying Lemma \ref{lem:X&DeltaX} with $X=\Lambda$ the action of $\Delta(\Lambda)$ on the $\U\otimes \U$-module $\C^{2^\Omega}(x_0)$ is identical to the action of $\Lambda$ on the $\U$-module $\C^{2^\Omega}$. Combined with Lemma \ref{lem:Lambda_2Omega} and (\ref{H:B}) this yields (\ref{B_2Omega(x0)}). The result follows.
\end{proof}

For the rest of this paper we always assume that $D\geq 2$. Suppose that $k$ is an integer with $1\leq k\leq D-1$. 
The {\it Johnson graph} $J(D,k)$ is a finite simple connected graph whose vertex set is ${\Omega \choose k}$ and two vertices $x,y$ are adjacent whenever $x\cap y\subsetdot x$. Hence the adjacency operator $\A$ of $J(D,k)$ is a linear endomorphism of $\C^{{\Omega \choose k}}$ given by 
\begin{gather}\label{ad:Johnson}
\A x
=
\sum_{\substack{|x|=|y|\\ x\cap y\subsetdot x}}
y 
\qquad 
\hbox{for all $x\in {\Omega \choose k}$}.
\end{gather}
Suppose that $x_0\in {\Omega \choose k}$. By \cite{BannaiIto1984,TerAlgebraI,TerAlgebraII,TerAlgebraIII} the dual adjacency operator $\A^*(x_0)$ of $J(D,k)$ with respect to $x_0$ is a linear endomorphism of $\C^{{\Omega \choose k}}$ given by 
\begin{gather}\label{dual ad:Johnson}
\A^*(x_0) x 
=
(D-1)
\left(
1-\frac{D(|x_0\setminus x|+|x\setminus x_0|)}{2k(D-k)}
\right) x
\qquad 
\hbox{for all $x\in {\Omega \choose k}$}.
\end{gather}
The {\it Terwilliger algebra $\T(x_0)$ of $J(D,k)$ with respect to $x_0$} is the subalgebra of ${\rm End}(\C^{{\Omega \choose k}})$ generated by $\A$ and $\A^*(x_0)$.
The algebra $\T(x_0)$ is a finite-dimensional semisimple algebra. Since $\T(x_0)\subseteq {\rm End}(\C^{{\Omega\choose k}})$ every irreducible $\T(x_0)$-module is contained in the $\T(x_0)$-module $\C^{{\Omega\choose k}}$ up to isomorphism. 
Therefore the algebra $\T(x_0)$ is isomorphic to
\begin{gather}\label{decT}
\bigoplus{\rm End}(\C^{\dim V})
\end{gather}
where the direct sum is over all non-isomorphic irreducible $\T(x_0)$-submodules $V$ of $\C^{{\Omega\choose k}}$. 
Please see \cite{TerAlgebraI,TerAlgebraII,TerAlgebraIII,Wedderburn1962} for details.

\begin{thm}
\label{thm:H->T}
Suppose that $k$ is an integer with $1\leq k\leq D-1$. 
For any $x_0\in {\Omega\choose k}$ the following equation holds:
$$
\T(x_0)={\rm Im}
\left(
\H\to {\rm End}(\C^{{\Omega \choose k}})
\right),
$$
where $\H\to {\rm End}(\C^{{\Omega \choose k}})$ is the representation of $\H$ into ${\rm End}(\C^{{\Omega \choose k}})$ corresponding to the $\H$-module $\C^{{\Omega \choose k}}(x_0)$.
\end{thm}
\begin{proof}
Comparing Theorem \ref{thm:H&2Omega} with (\ref{ad:Johnson}) and (\ref{dual ad:Johnson}), the following equations hold on the $\H$-module $\C^{{\Omega \choose k}}(x_0)$:
\begin{align*}
\A&=B-\frac{D}{2}-\frac{(D-2k)^2}{4},
\\
\A^*(x_0) &= 
\frac{D(D-1)}{k(D-k)}
\left(
A-\frac{(D-2k)^2}{4D}
\right).
\end{align*}
The result follows.
\end{proof}


\section{The dimensions of the Terwilliger algebras of Johnson graphs}\label{s:dim}

Observe that the graph $J(D,k)$ is isomorphic to $J(D,D-k)$. Without loss of generality we may assume that $1\leq k\leq \frac{D}{2}$. Let $x_0\in {\Omega\choose k}$ be given. 
From the viewpoint (\ref{decT}) the classification of irreducible $\T(x_0)$-modules up to isomorphism plays a crucial role in the study of $\T(x_0)$. 
It was stated in \cite[Example 6.1]{TerAlgebraIII} without proof that the non-isomorphic irreducible $\T(x_0)$-modules are classified by three parameters called the endpoints, the dual endpoints, and the diameters. 
Let 
\begin{gather}
\label{I(k)}
\I(k)=\left\{(i,j)\in \Z^2
\,\Big|\, 0\leq i\leq \frac{k}{2},0\leq j\leq \min\left\{k,\frac{D-k}{2}\right\}, 0\leq i+j\leq k
\right\}.
\end{gather}
A few years ago it was noticed in \cite{Johnson:2017} that the parametrization of the isomorphism classes of irreducible $\T(x_0)$-modules can be simplified to the set $\I(k)$ 
when $k<\frac{D}{2}$; it can be simplified to $\{(i,j)\in \I(\frac{D}{2})\,|\,j\leq i\}=\{(i,j)\in\Z^2\,|\, 0\leq j\leq i\leq \frac{D}{4}\}$ when $k=\frac{D}{2}$. 
By giving a character theoretical explanation for $\I(k)$, the main result of \cite{Johnson:2019} proved that $\T(x_0)$ is equal to the centralizer algebra for the stabilizer of $x_0$ in the automorphism group of $J(D,k)$.

For any integer $i$ with $0\leq i\leq k$  the $i^{\rm\, th}$ dual primitive idempotent $\mathbf E_i^*(x_0)$ of $J(D,k)$ with respect to $x_0$ is a linear endomorphism of $\C^{{\Omega\choose k}}$ defined by 
$$
\mathbf E_i^*(x_0)x=
\left\{
\begin{array}{ll}
x
\qquad &\hbox{if $|x\cap x_0|=k-i$},
\\
0
\qquad &\hbox{else}
\end{array}
\right.
\qquad 
\hbox{for all $x\in {\Omega\choose k}$}.
$$
In \cite{Johnson:2014,Johnson:2007} the algebra $\T(x_0)$ was looked into from the perspective of the decomposition:
$$
\T(x_0)=\bigoplus_{i=0}^k \bigoplus_{j=0}^k \mathbf E_i^*(x_0) \T(x_0)\mathbf E_j^*(x_0).
$$
After complicated analysis, it was claimed in  \cite{Johnson:2014,Johnson:2007} that the structure of the algebra $\T(x_0)$ is  determined.
However the correctness was denied in \cite{Johnson:2017}. The correct result is as follows:

\begin{prop}
[page 152,\cite{Johnson:2017}]
\label{prop:Ito}
Suppose that $k$ is an integer with $1\leq k\leq \frac{D}{2}$. Then the following statements hold:
\begin{enumerate}
\item Suppose that $1\leq k<\frac{D}{2}$. For any $x_0\in {\Omega\choose k}$ the algebra $\T(x_0)$ is isomorphic to 
\begin{align*}
\bigoplus_{(i,j)\in \I_{\rm I}(k)}
{\rm End}
(
\C^{k-2i+1}
)
\oplus 
\bigoplus_{(i,j)\in \I_{\rm II}(k)}
{\rm End}
(
\C^{k-i-j+1}
)
\oplus 
\bigoplus_{(i,j)\in \I_{\rm III}(k)}
{\rm End}
(
\C^{D-k-2j+1}
)
\end{align*}
where 
\begin{align*}
\I_{\rm I}(k)
&=\left\{ 
(i,j)\in \I(k)
\,|\,
j\leq i
\right\},
\\
\I_{\rm II}(k)
&=\left\{ 
(i,j)\in \I(k)
\,|\,
0< j-i\leq D-2k
\right\},
\\
\I_{\rm III}(k)
&=\left\{ 
(i,j)\in \I(k)
\,|\,
D-2k< j-i
\right\}.
\end{align*}

\item Suppose that $k=\frac{D}{2}$. For any $x_0\in {\Omega\choose k}$ the algebra $\T(x_0)$ is isomorphic to 
\begin{align*}
\bigoplus_{i=0}^{\floor{\frac{D}{4}}}
\,
(i+1)
\cdot 
{\rm End}
(
\C^{\frac{D}{2}-2i+1}
).
\end{align*}
\end{enumerate}
\end{prop}

In the final section we shall study the structure of $\T(x_0)$ along our own veins and express the dimension of $\T(x_0)$ in terms of binomials coefficients.
To do this we will use the following common equations:
\begin{align}
\sum\limits_{i=0}^{\floor{\frac{n}{2}}} (n-2i)^2&={n+2\choose 3}
\qquad 
\hbox{for all $n\in \N$},
\label{b1}
\\
\sum\limits_{i=0}^n (n-i)^2&={n+2\choose 3}+{n+1\choose 3}
\qquad 
\hbox{for all $n\in \N$},
\label{b0}
\\
\sum\limits_{i=0}^n {i\choose \ell}&={n+1\choose \ell+1}
\qquad 
\hbox{for all $\ell,n\in \N$}.
\label{b2}
\end{align}
Define 
\begin{gather}\label{sln}
s_{\ell}(n)=\sum_{i =0}^{\floor{\frac{n}{2}}}{n-2i \choose \ell}
\qquad 
\hbox{for all $\ell,n\in \N$}.
\end{gather}

\begin{lem}\label{lem:sln}
For all $\ell,n\in \N$ the following equations hold:
\begin{enumerate}
\item $s_{\ell+1}(n+1)+s_{\ell+1}(n)={n+2\choose \ell+2}$.

\item $s_{\ell+1}(n+1)-s_{\ell+1}(n)=s_{\ell}(n)$.
\end{enumerate}
\end{lem} 
\begin{proof}
(i): By (\ref{sln}) the left-hand side of Lemma \ref{lem:sln}(i) is equal to $\sum\limits_{i=0}^{n+1}{i\choose \ell+1}$. By  (\ref{b2}) it is equal to the right-hand side of Lemma \ref{lem:sln}(i).

(ii): Immediate from Pascal's rule.
\end{proof}

Recall that for any real number $a$, the notation $\ceil{a}$ stands for the least integer greater than or equal to $a$ and the notation $\floor{a}$ stands for the greatest integer less than or equal to $a$. 

\begin{lem}\label{lem1:sln}
For all $\ell,n\in \N$ the following equations hold:
\begin{enumerate}
\item 
\begin{gather*}
\begin{split}
s_{\ell}(n)=
\frac{1}{2}{n+1\choose \ell+1}
+
\frac{1}{4}
\sum\limits_{i=0}^{\ell-1}
\left(\frac{-1}{2}\right)^i
{n+1\choose \ell-i}
+
\left\{
\begin{array}{ll}
\frac{(-1)^\ell}{2^{\ell+1}}
\qquad 
&\hbox{if $n$ is even},
\\
0
\qquad 
&\hbox{if $n$ is odd}.
\end{array}
\right.
\end{split}
\end{gather*}

\item 
\begin{gather*}
\begin{split}
s_{\ell}(n)=
\left\{
\begin{array}{ll}
\displaystyle
\ceil{
\frac{1}{2}{n+1\choose \ell+1}
+
\frac{1}{4}
\sum\limits_{i=0}^{\ell-1}
\left(\frac{-1}{2}\right)^i
{n+1\choose \ell-i}
}
\qquad 
&\hbox{if $\ell$ is even},
\\
\displaystyle
\floor{
\frac{1}{2}{n+1\choose \ell+1}
+
\frac{1}{4}
\sum\limits_{i=0}^{\ell-1}
\left(\frac{-1}{2}\right)^i
{n+1\choose \ell-i}
}
\qquad 
&\hbox{if $\ell$ is odd}.
\end{array}
\right.
\end{split}
\end{gather*}
\end{enumerate}
\end{lem}
\begin{proof}
(i): Let $n\in \N$ be given. Subtracting Lemma \ref{lem:sln}(ii) from Lemma \ref{lem:sln}(i) yields that 
\begin{gather}\label{s3}
2s_{\ell+1}(n)+s_{\ell}(n)={n+2\choose \ell+2}
\qquad 
\hbox{for all $\ell\in \N$}.
\end{gather}
Hence the sequence $\{s_{\ell}(n)\}_{\ell\in \N}$ is uniquely determined by the recurrence (\ref{s3}) with the initial condition $s_{0}(n)=\floor{\frac{n}{2}}+1$. Let $s'_{\ell}(n)$ denote the right-hand side of Lemma \ref{lem1:sln}(i) for all $\ell\in \N$. Observe that 
$$
s'_0(n)
=\left\{
\begin{array}{ll}
\frac{n}{2}+1
\qquad 
&\hbox{if $n$ is even},
\\
\frac{n+1}{2}
\qquad 
&\hbox{if $n$ is odd}.
\end{array}
\right.
$$
Hence $s'_0(n)=s_0(n)$. A direct calculation shows that the sequence $\{s'_{\ell}(n)\}_{\ell\in \N}$ satisfies the recurrence (\ref{s3}).
Hence $s'_{\ell}(n)=s_{\ell}(n)$ for all $\ell\in \N$ and Lemma \ref{lem1:sln}(i) follows. 

(ii): Immediate from Lemma \ref{lem1:sln}(i).
\end{proof}

Recall the index set $\P(k)$ from (\ref{P}). 
Observe that 
\begin{gather}\label{P:k<D/2}
\P(k)=
\left\{
(i,j)\in \Z^2
\,\Big|\,
0\leq i\leq \frac{D-k}{2}, 
0\leq j\leq \min\left\{k-i,\frac{k}{2}\right\}
\right\}
\qquad 
\hbox{if $k<\displaystyle \frac{D}{2}$}.
\end{gather}

\begin{prop}\label{prop1:T}
Suppose that $k$ is an integer with $1\leq k< \frac{D}{3}$. For any $x_0\in {\Omega\choose k}$ the following statements hold: 
\begin{enumerate}
\item The algebra $\T(x_0)$ is isomorphic to 
\begin{align*}
\bigoplus_{(i,j)\in \P_{{\rm I}}(k)}
{\rm End}(\C^{k-2j+1})
\oplus
\bigoplus_{(i,j)\in \P_{{\rm II}}(k)}
{\rm End}(\C^{k-i-j+1})
\end{align*}
where 
\begin{align*}
\P_{{\rm I}}(k)
&=
\left\{
(i,j)\in \Z^2
\,\Big|\,
0\leq i\leq j\leq \frac{k}{2}
\right\},
\\
\P_{{\rm II}}(k)
&=
\left\{
(i,j)\in \Z^2
\,\Big|\,
0\leq j<i\leq \frac{k}{2}
\right\}
\cup
\left\{
(i,j)\in \Z^2
\,\Big|\,
\frac{k}{2}<i\leq k,0\leq j\leq k-i
\right\}.
\end{align*}

\item The dimension of $\T(x_0)$ is equal to 
$$
{k+3\choose 4}
+
\floor{\frac{1}{2}{k+4\choose 4}
+\frac{1}{4}{k+4\choose 3}
-\frac{1}{8}{k+4\choose 2}
+\frac{1}{16}{k+4\choose 1}}.
$$
\end{enumerate}
\end{prop}
\begin{proof}
(i): 
Using (\ref{P:k<D/2})
it is routine to verify that 
\begin{align*}
\P_{{\rm I}}(k)&=\{(i,j)\in \P(k)\,|\, i\leq j\},
\\
\P_{{\rm II}}(k)&=\{(i,j)\in \P(k)\,|\,j<i\}.
\end{align*}
Let $(i,j)\in \P_{{\rm I}}(k)$ be given.
Applying Theorem \ref{thm:Lmn}(i) with $(m,n,\ell)=(D-k-2i,k-2j,k-i-j)$ yields that 
$\dim\, (L_{D-k-2i}\otimes L_{k-2j})(D-2k)=k-2j+1$. 
Let $(i,j)\in \P_{{\rm II}}(k)$ be given. Applying Theorem \ref{thm:Lmn}(i) with $(m,n,\ell)=(D-k-2i,k-2j,k-i-j)$ yields that 
$\dim (L_{D-k-2i}\otimes L_{k-2j})(D-2k)=k-i-j+1$. 
Combined with Theorems \ref{thm:H&Omega k}(i), \ref{thm:H->T} and (\ref{decT}) the statement (i) follows.

(ii): By Proposition \ref{prop1:T}(i) the dimension of $\T(x_0)$ is equal to 
\begin{gather}\label{T1}
\sum_{i=0}^{\floor{\frac{k}{2}}}
\sum_{j=i}^{\floor{\frac{k}{2}}}
(k-2j+1)^2
+
\sum_{j=0}^{\floor{\frac{k}{2}}}
\sum_{i=j+1}^{\floor{\frac{k}{2}}}
(k-i-j+1)^2
+
\sum_{i=\floor{\frac{k}{2}}+1}^k
\sum_{j=0}^{k-i}
(k-i-j+1)^2.
\end{gather}
Using (\ref{b1}) and (\ref{sln}) yields that the first summation in (\ref{T1}) is equal to $
\sum\limits_{i=0}^{\floor{\frac{k}{2}}}
{k-2i+3\choose 3}=s_3(k+3)$.
Using (\ref{b0})--(\ref{sln}) yields that the second summation in (\ref{T1}) is equal  to 
\begin{eqnarray*}
&&\sum_{j=0}^{\floor{\frac{k}{2}}}
{k-2j+2\choose 3}
+
{k-2j+1\choose 3}
-
{\ceil{\frac{k}{2}}-j+2\choose 3}
-
{\ceil{\frac{k}{2}}-j+1\choose 3}
\\
&=&
s_3(k+2)
+s_3(k+1)
-{\ceil{\frac{k}{2}}+2\choose 4}
-{\ceil{\frac{k}{2}}+3\choose 4}.
\end{eqnarray*}
Using (\ref{b0}) and (\ref{b2}) yields that the third summation in (\ref{T1}) is equal  to 
\begin{eqnarray*}
\sum_{i=\floor{\frac{k}{2}}+1}^{k}
{k-i+3\choose 3}
+
{k-i+2\choose 3}
=
{\ceil{\frac{k}{2}}+3\choose 4}
+{\ceil{\frac{k}{2}}+2\choose 4}.
\end{eqnarray*}
Hence $\dim \T(x_0)=s_3(k+3)+s_3(k+2)+s_3(k+1)$. By Lemma \ref{lem:sln}(i) the sum $s_3(k+2)+s_3(k+1)={k+3\choose 4}$. 
Combined with Lemma \ref{lem1:sln}(ii) the statement (ii) follows.
\end{proof}

\begin{prop}\label{prop2:T}
Suppose that $k$ is an integer with $\frac{D}{3}\leq k< \frac{2D}{5}$. For any $x_0\in {\Omega\choose k}$ the following statements hold: 
\begin{enumerate}
\item The algebra $\T(x_0)$ is isomorphic to 
\begin{align*}
\bigoplus_{(i,j)\in \P_{{\rm I}}(k)}
{\rm End}(\C^{k-2j+1})
\oplus
\bigoplus_{(i,j)\in \P_{{\rm II}}(k)}
{\rm End}(\C^{k-i-j+1})
\oplus
\bigoplus_{(i,j)\in \P_{{\rm III}}(k)}
{\rm End}(\C^{D-k-2i+1})
\end{align*}
where 
\begin{align*}
\P_{{\rm I}}(k)
&=
\left\{
(i,j)\in \Z^2\,
\,\Big|\,
0\leq i\leq j\leq \frac{k}{2}
\right\},
\\
\P_{{\rm II}}(k)
&=
\left\{
(i,j)\in \Z^2
\,\Big|\,
0\leq j<i\leq \frac{k}{2}
\right\}
\\
&\qquad \cup  
\left\{
(i,j)\in \Z^2
\,\Big|\,
\frac{k}{2}<i\leq D-2k,0\leq j\leq k-i
\right\}
\\
&\qquad \cup 
\left\{(i,j)\in \Z^2
\,\Big|\,
D-2k<i\leq \frac{D-k}{2},2k-D+i\leq j\leq k-i
\right\},
\\
\P_{{\rm III}}(k)
&=
\left\{
(i,j)\in \Z^2
\,\Big|\,
0\leq j< 2k-D+i\leq \frac{3k-D}{2}
\right\}.
\end{align*}

\item The dimension of $\T(x_0)$ is equal to 
\begin{align*}
&{k+3\choose 4}
+\floor{\frac{1}{2}{k+4\choose 4}
+\frac{1}{4}{k+4\choose 3}
-\frac{1}{8}{k+4\choose 2}
+\frac{1}{16}{k+4\choose 1}}
\\
&\quad -
\floor{\frac{1}{2}{3k-D+3\choose 4}
+\frac{1}{4}{3k-D+3\choose 3}
-\frac{1}{8}{3k-D+3\choose 2}
+\frac{1}{16}{3k-D+3\choose 1}}.
\end{align*}
\end{enumerate}
\end{prop}
\begin{proof}
(i): 
Using (\ref{P:k<D/2})
it is routine to verify that 
\begin{align*}
\P_{{\rm I}}(k)&=\{(i,j)\in \P(k)\,|\, i\leq j\},
\\
\P_{{\rm II}}(k)&=\{(i,j)\in \P(k)\,|\, 2k-D+i\leq j<i\},
\\
\P_{{\rm III}}(k)&=\{(i,j)\in \P(k)\,|\, j<2k-D+i\}.
\end{align*}
Let $(i,j)\in \P_{{\rm I}}(k)$ be given.  
Applying Theorem \ref{thm:Lmn}(i) with $(m,n,\ell)=(D-k-2i,k-2j,k-i-j)$ yields that $\dim (L_{D-k-2i}\otimes L_{k-2j})(D-2k)=k-2j+1$. 
Let $(i,j)\in \P_{{\rm II}}(k)$ be given.  
Applying Theorem \ref{thm:Lmn}(i) with $(m,n,\ell)=(D-k-2i,k-2j,k-i-j)$ yields that $\dim (L_{D-k-2i}\otimes L_{k-2j})(D-2k)=k-i-j+1$. 
Let $(i,j)\in \P_{{\rm III}}(k)$ be given.  
Applying Theorem \ref{thm:Lmn}(i) with $(m,n,\ell)=(D-k-2i,k-2j,k-i-j)$ yields that $\dim (L_{D-k-2i}\otimes L_{k-2j})(D-2k)=D-k-2i+1$. 
Combined with Theorems \ref{thm:H&Omega k}(i), \ref{thm:H->T} and (\ref{decT}) the statement (i) follows.

(ii): By Proposition \ref{prop2:T}(i) the dimension of $\T(x_0)$ is equal to 
\begin{align}\label{T2}
\begin{split}
&\sum_{i=0}^{\floor{\frac{k}{2}}}
\sum_{j=i}^{\floor{\frac{k}{2}}}
(k-2j+1)^2
+
\sum_{j=0}^{\floor{\frac{k}{2}}}
\sum_{i=j+1}^{\floor{\frac{k}{2}}}
(k-i-j+1)^2
+
\sum_{i=\floor{\frac{k}{2}}+1}^{D-2k}
\sum_{j=0}^{k-i}
(k-i-j+1)^2
\\
&\qquad +
\sum_{i=D-2k+1}^{\floor{\frac{D-k}{2}}}
\sum_{j=2k-D+i}^{k-i}
(k-i-j+1)^2
+
\sum_{j=0}^{\floor{\frac{3k-D}{2}}}
\sum_{i=D-2k+j+1}^{\floor{\frac{D-k}{2}}}
(D-k-2i+1)^2.
\end{split}
\end{align}
Using (\ref{b1}) and (\ref{sln}) yields that the first summation in (\ref{T2}) is equal to $s_3(k+3)$.
Using (\ref{b0})--(\ref{sln}) yields that the second summation in (\ref{T2}) is equal  to 
$
s_3(k+1)
+s_3(k+2)
-{\ceil{\frac{k}{2}}+2\choose 4}
-{\ceil{\frac{k}{2}}+3\choose 4}$.
Using (\ref{b0}) and (\ref{b2}) yields that the third summation in (\ref{T2}) is equal to ${\ceil{\frac{k}{2}}+2\choose 4}+{\ceil{\frac{k}{2}}+3\choose 4}-{3k-D+3\choose 4}-{3k-D+2\choose 4}$. Using (\ref{b0}) and (\ref{sln}) yields that the fourth summation in (\ref{T2}) is equal to $s_3(3k-D)+s_3(3k-D+1)$. 
Using (\ref{b1}) and (\ref{sln}) yields the fifth summation in (\ref{T2}) is equal to $s_3(3k-D+1)$. Combined with Lemma \ref{lem:sln}(i) the dimension of $\T(x_0)$ is equal to 
$$
{k+3\choose 4}+s_3(k+3)
-s_3(3k-D+2).
$$
Combined with Lemma \ref{lem1:sln}(ii) the statement (ii) follows.
\end{proof}

\begin{prop}\label{prop3:T}
Suppose that $k$ is an integer with $\frac{2D}{5}\leq k< \frac{D}{2}$. For any $x_0\in {\Omega\choose k}$ the following statements hold: 
\begin{enumerate}
\item The algebra $\T(x_0)$ is isomorphic to 
\begin{align*}
\bigoplus_{(i,j)\in \P_{{\rm I}}(k)}
{\rm End}(\C^{k-2j+1})
\oplus
\bigoplus_{(i,j)\in \P_{{\rm II}}(k)}
{\rm End}(\C^{k-i-j+1})
\oplus
\bigoplus_{(i,j)\in \P_{{\rm III}}(k)}
{\rm End}(\C^{D-k-2i+1})
\end{align*}
where 
\begin{align*}
\P_{{\rm I}}(k)
&=
\left\{
(i,j)\in \Z^2
\,\Big|\,
0\leq i\leq j\leq \frac{k}{2}
\right\},
\\
\P_{{\rm II}}(k)
&=
\left\{
(i,j)\in \Z^2
\,\Big|\,
0\leq j<i\leq D-2k
\right\}
\\
&\qquad \cup 
\left\{
(i,j)\in \Z^2
\,\Big|\,
D-2k< i\leq \frac{k}{2}, 2k-D+i\leq j< i
\right\}
\\
&\qquad \cup 
\left
\{(i,j)\in \Z^2
\,\Big|\,
\frac{k}{2}< i\leq \frac{D-k}{2}, 2k-D+i\leq j\leq k-i
\right\},
\\
\P_{{\rm III}}(k)
&=
\left\{
(i,j)\in \Z^2
\,\Big|\,
0\leq j< 2k-D+i\leq \frac{3k-D}{2}
\right\}.
\end{align*}

\item The dimension of $\T(x_0)$ is equal to 
\begin{align*}
&{k+3\choose 4}
+\floor{\frac{1}{2}{k+4\choose 4}
+\frac{1}{4}{k+4\choose 3}
-\frac{1}{8}{k+4\choose 2}
+\frac{1}{16}{k+4\choose 1}}
\\
&\quad -
\floor{\frac{1}{2}{3k-D+3\choose 4}
+\frac{1}{4}{3k-D+3\choose 3}
-\frac{1}{8}{3k-D+3\choose 2}
+\frac{1}{16}{3k-D+3\choose 1}}.
\end{align*}
\end{enumerate}
\end{prop}
\begin{proof}
(i): Using (\ref{P:k<D/2})
it is routine to verify that 
\begin{align*}
\P_{{\rm I}}(k)&=\{(i,j)\in \P(k)\,|\, i\leq j\},
\\
\P_{{\rm II}}(k)&=\{(i,j)\in \P(k)\,|\, 2k-D+i\leq j<i\},
\\
\P_{{\rm III}}(k)&=\{(i,j)\in \P(k)\,|\, j<2k-D+i\}.
\end{align*}
Let $(i,j)\in \P_{{\rm I}}(k)$ be given.  
Applying Theorem \ref{thm:Lmn}(i) with $(m,n,\ell)=(D-k-2i,k-2j,k-i-j)$ yields that $\dim (L_{D-k-2i}\otimes L_{k-2j})(D-2k)=k-2j+1$. 
Let $(i,j)\in \P_{{\rm II}}(k)$ be given.  
Applying Theorem \ref{thm:Lmn}(i) with $(m,n,\ell)=(D-k-2i,k-2j,k-i-j)$ yields that $\dim (L_{D-k-2i}\otimes L_{k-2j})(D-2k)=k-i-j+1$. 
Let $(i,j)\in \P_{{\rm III}}(k)$ be given.  
Applying Theorem \ref{thm:Lmn}(i) with $(m,n,\ell)=(D-k-2i,k-2j,k-i-j)$ yields that $\dim (L_{D-k-2i}\otimes L_{k-2j})(D-2k)=D-k-2i+1$. 
Combined with Theorems \ref{thm:H&Omega k}(i), \ref{thm:H->T} and (\ref{decT}) the statement (i) follows.

(ii): By Proposition \ref{prop3:T}(i) the dimension of $\T(x_0)$ is equal to 
\begin{align}\label{T3}
\begin{split}
&\sum_{i=0}^{\floor{\frac{k}{2}}}
\sum_{j=i}^{\floor{\frac{k}{2}}}
(k-2j+1)^2
+
\sum_{j=0}^{D-2k}
\sum_{i=j+1}^{D-2k}
(k-i-j+1)^2
+
\sum_{i=D-2k+1}^{\floor{\frac{k}{2}}}
\sum_{j=2k-D+i}^{i-1}
(k-i-j+1)^2
\\
&\qquad +
\sum_{i=\floor{\frac{k}{2}}+1}^{\floor{\frac{D-k}{2}}}
\sum_{j=2k-D+i}^{k-i}
(k-i-j+1)^2
+
\sum_{j=0}^{\floor{\frac{3k-D}{2}}}
\sum_{i=D-2k+j+1}^{\floor{\frac{D-k}{2}}}
(D-k-2i+1)^2.
\end{split}
\end{align}
Using (\ref{b1}) and (\ref{sln}) yields that the first summation in (\ref{T3}) is equal to $s_3(k+3)$.
Using (\ref{b0})--(\ref{sln}) yields that the second summation in (\ref{T3}) is equal  to 
$
s_3(k+1)
+s_3(k+2)
-s_3(5k-2D-1)
-s_3(5k-2D)
+{5k-2D+1\choose 4}
+{5k-2D+2\choose 4}
-{3k-D+3\choose 4}
-{3k-D+2\choose 4}$.
Using (\ref{b0}) and (\ref{sln}) yields that the third summation in (\ref{T3}) is equal to 
$
s_3(3k-D)
+s_3(3k-D+1)
-s_3(5k-2D)
-s_3(5k-2D+1)
-s_3(D-k-2\floor{\frac{k}{2}})
-s_3(D-k-2\floor{\frac{k}{2}}+1)
$
and the fourth summation in (\ref{T3}) is equal to 
$
s_3(D-k-2\floor{\frac{k}{2}})
+s_3(D-k-2\floor{\frac{k}{2}}+1)$. 
Using (\ref{b1}) and (\ref{sln}) yields the fifth summation in (\ref{T3}) is equal to $s_3(3k-D+1)$. Combined with Lemma \ref{lem:sln}(i) the dimension of $\T(x_0)$ is equal to 
$$
{k+3\choose 4}+s_3(k+3)-s_3(3k-D+2).
$$
Combined with Lemma \ref{lem1:sln}(ii) the statement (ii) follows.
\end{proof}

\begin{prop}\label{prop4:T}
Suppose that $k=\frac{D}{2}$ is an integer. For any $x_0\in {\Omega\choose k}$ the following statements hold: 
\begin{enumerate}
\item The algebra $\T(x_0)$ is isomorphic to 
\begin{align}
\label{T_D/2}
\bigoplus_{i=0}^{\floor{\frac{D}{4}}}
\bigoplus_{j=i}^{\floor{\frac{D}{4}}}
{\rm End}(\C^{\frac{D}{2}-2j+1}).
\end{align}

\item The dimension of $\T(x_0)$ is equal to 
\begin{align*}
\floor{\frac{1}{2}{\frac{D}{2}+4\choose 4}
+\frac{1}{4}{\frac{D}{2}+4\choose 3}
-\frac{1}{8}{\frac{D}{2}+4\choose 2}
+\frac{1}{16}{\frac{D}{2}+4\choose 1}}.
\end{align*}
\end{enumerate}
\end{prop}
\begin{proof}
(i): Immediate from Theorems \ref{thm:Lmn}(i), \ref{thm:H&Omega k}(ii), \ref{thm:H->T} and (\ref{decT}).

(ii): By Proposition \ref{prop4:T}(i) the dimension of $\T(x_0)$ is equal to 
\begin{gather*}
\sum_{i=0}^{\floor{\frac{D}{4}}}
\sum_{j=i}^{\floor{\frac{D}{4}}}
\left(
\frac{D}{2}-2j+1
\right)^2.
\end{gather*}
Using (\ref{b1}) and (\ref{sln}) yields the double summation is equal to $s_3(\frac{D}{2}+3)$. Combined with Lemma \ref{lem1:sln}(ii) the statement (ii) follows.
\end{proof}

\begin{rem}
(i): Suppose that $k$ is an integer with $1\leq k<\frac{D}{2}$.
By (\ref{I(k)}) and (\ref{P:k<D/2}) there is a bijection $\I(k)\to \P(k)$ given by 
\begin{eqnarray*}
(i,j) &\mapsto & (j,i)
\qquad 
\hbox{for all $(i,j)\in \I(k)$}.
\end{eqnarray*}
Recall the sets $\I_{\rm I}(k), \I_{\rm II}(k), \I_{\rm III}(k)$ from Proposition \ref{prop:Ito}(i) and the sets $\P_{\rm I}(k), \P_{\rm II}(k), \P_{\rm III}(k)$ from Propositions \ref{prop1:T}(i)--\ref{prop3:T}(i). Note that $\I_{\rm III}(k)=\emptyset$ when $k<\frac{D}{3}$. From the proofs of Propositions \ref{prop1:T}--\ref{prop3:T}, it follows that the bijection $\I(k)\to \P(k)$ maps $\I_{\rm I}(k), \I_{\rm II}(k), \I_{\rm III}(k)$ to $\P_{\rm I}(k), \P_{\rm II}(k), \P_{\rm III}(k)$ respectively. Therefore Propositions \ref{prop1:T}(i)--\ref{prop3:T}(i) coincide with Proposition \ref{prop:Ito}(i).

(ii): For any integer $i$ with $0\leq i\leq \floor{\frac{D}{4}}$ the number of ${\rm End}(\C^{\frac{D}{2}-2i+1})$ in (\ref{T_D/2}) is equal to $i+1$. Therefore Proposition \ref{prop4:T}(i) coincides with Proposition \ref{prop:Ito}(ii).
\end{rem}

\begin{thm}\label{thm:dimT}
Suppose that $k$ is an integer with $1\leq k\leq D-1$. 
For any $x_0\in {\Omega\choose k}$ the following statements hold:
\begin{enumerate}
\item Suppose that $1\leq k< \frac{D}{3}$. 
Then the dimension of $\T(x_0)$ is equal to 
\begin{align*}
{k+3\choose 4}
+
\floor{\frac{1}{2}{k+4\choose 4}
+\frac{1}{4}{k+4\choose 3}
-\frac{1}{8}{k+4\choose 2}
+\frac{1}{16}{k+4\choose 1}}.
\end{align*}

\item Suppose that $\frac{D}{3}\leq k<\frac{D}{2}$. 
Then the dimension of $\T(x_0)$ is equal to 
\begin{align*}
&{k+3\choose 4}
+\floor{\frac{1}{2}{k+4\choose 4}
+\frac{1}{4}{k+4\choose 3}
-\frac{1}{8}{k+4\choose 2}
+\frac{1}{16}{k+4\choose 1}}
\\
&\quad -
\floor{\frac{1}{2}{3k-D+3\choose 4}
+\frac{1}{4}{3k-D+3\choose 3}
-\frac{1}{8}{3k-D+3\choose 2}
+\frac{1}{16}{3k-D+3\choose 1}}.
\end{align*}

\item Suppose that $k=\frac{D}{2}$. 
Then the dimension of $\T(x_0)$ is equal to 
\begin{align*}
\floor{\frac{1}{2}{\frac{D}{2}+4\choose 4}
+\frac{1}{4}{\frac{D}{2}+4\choose 3}
-\frac{1}{8}{\frac{D}{2}+4\choose 2}
+\frac{1}{16}{\frac{D}{2}+4\choose 1}}.
\end{align*}

\item Suppose that $\frac{D}{2}<k\leq \frac{2D}{3}$. 
Then the dimension of $\T(x_0)$ is equal to 
\begin{align*}
&{D-k+3\choose 4}
+\floor{\frac{1}{2}{D-k+4\choose 4}
+\frac{1}{4}{D-k+4\choose 3}
-\frac{1}{8}{D-k+4\choose 2}
+\frac{1}{16}{D-k+4\choose 1}}
\\
&\quad -
\floor{\frac{1}{2}{2D-3k+3\choose 4}
+\frac{1}{4}{2D-3k+3\choose 3}
-\frac{1}{8}{2D-3k+3\choose 2}
+\frac{1}{16}{2D-3k+3\choose 1}}.
\end{align*}

\item Suppose that $\frac{2D}{3}< k\leq D-1$. 
Then the dimension of $\T(x_0)$ is equal to 
\begin{align*}
{D-k+3\choose 4}
+
\floor{\frac{1}{2}{D-k+4\choose 4}
+\frac{1}{4}{D-k+4\choose 3}
-\frac{1}{8}{D-k+4\choose 2}
+\frac{1}{16}{D-k+4\choose 1}}.
\end{align*}
\end{enumerate}
\end{thm}
\begin{proof}
(i)--(iii): Immediate from Propositions \ref{prop1:T}(ii)--\ref{prop4:T}(ii).

(iv), (v): Since the graph $J(D,k)$ is isomorphic to $J(D,D-k)$ the statements (iv) and (v) follow by replacing $k$ by $D-k$ in Theorem \ref{thm:dimT}(i), (ii).
\end{proof}

We conclude the paper with the relationship to the literature \cite{Johnson:2021,Watanabe:2017,Watanabe2020,Huang:CG}. 
The algebra $U_q(\mathfrak{sl}_2)$ is a $q$-analog of $\U$. 
The algebraic treatment of the Clebsch--Gordan coefficients of $U_q(\mathfrak{sl}_2)$ was proposed in \cite[Theorem 2.9]{Huang:CG} earlier than Theorem \ref{thm:H->U2}. 
In \cite{Johnson:2021}, applying the theory of angular momentum coupling yielded a decomposition of $\T(x_0)$-module $\C^{{\Omega\choose k}}$ into a direct sum of irreducible $\T(x_0)$-modules over a vague index set and the classification of irreducible $\T(x_0)$-modules up to isomorphism was ignored. 
Inspired by \cite{Huang:CG, Johnson:2021}, we first develop Theorems \ref{thm:H->U2} and \ref{thm:Lmn} to link the Clebsch--Gordan coefficients of $\U$ and the study of $\T(x_0)$ thoroughly. Theorem \ref{thm:dimT} is also first released in the present paper.

The Grassmann graph is a $q$-analog of the Johnson graph. The subspace lattice over a finite field is a $q$-analog of the subset lattice. 
In \cite{Watanabe2020,Watanabe:2017} it was mentioned how the Terwilliger algebra of the Grassmann graph is related to an algebra associated with a subspace lattice over a finite field. Comparing \cite{Huang:CG,Watanabe2020,Watanabe:2017} with the present paper, it still lacks some links between the Clebsch--Gordan coefficients of $U_q(\mathfrak{sl}_2)$ and the Terwilliger algebras of Grassmann graphs. The missing links will show up in the subsequent paper \cite{Huang:CG&Grassmann}.

\subsection*{Acknowledgements}
The research is supported by the National Science and Technology Council of Taiwan under the project NSTC 112-2115-M-008-009-MY2.

\bibliographystyle{amsplain}
\bibliography{MP}

\end{document}